\input amstex
\documentstyle{amsppt}
\topmatter
\magnification=\magstep1
\pagewidth{5.2 in}
\pageheight{6.7 in}
%\hcorrection{-0.4in}
%\vcorrection{-0.4in}
\abovedisplayskip=10pt
\belowdisplayskip=10pt
\NoBlackBoxes
\title
A note on the alternating sums of powers of consecutive integers
\endtitle
%Use \endgraf to indicate new paragraph
%\thanks will become a 1st page footnote
\author  Taekyun Kim \endauthor
\affil\rm{{Institute of Science Education,}\\
{ Kongju National University, Kongju 314-701, S. Korea}\\
{e-mail: tkim$\@$kongju.ac.kr}}\\
\endaffil

\abstract{ For $n,k\in\Bbb Z_{\geq 0},$ let $T_{n}(k)$ be the
alternating sums of the $n$-th powers of positive integers up to
$k-1$: $T_{n}(k)=\sum_{l=0}^{k-1}(-1)^ll^n .$ Following an idea
due to Euler, we give the below formula for $T_{n}(k)$:
$$ T_{n}(k)=\frac{(-1)^{k+1}}{2}\sum_{l=0}^{n-1}
\binom{n}{l}E_lk^{n-l}+\frac{E_n}{2}\left(1+(-1)^{k+1}\right),$$
where $E_l$ are the Euler numbers. }
\endabstract
\thanks 2000 Mathematics Subject Classification  11S80, 11B68 \endthanks
\thanks Key words and phrases: Euler number,zeta function, Bernoulli numbers \endthanks
\rightheadtext{ Taekyun Kim    } \leftheadtext{ Alternating Sums
of powers of consecutive integers }
%\TagsOnRight
\endtopmatter

\document

\head 1. Introduction \endhead
 J. Bernoulli (1713) first
discovered the method which one can produce those formulae for the
sum $\sum_{l=1}^nl^k$ for any natural numbers. The Bernoulli
numbers are among the most interesting and important number
sequences in mathematics. They first appeared in the posthumous
work ``ARS Conjectandi''(1713) by J. Bernoulli (1654-1705) in
connection with sums of powers of consecutive integers. Bernoulli
numbers are particularly important in number theory, especially in
connection with Fermat's last theorem. The number sequences of
Euler, Genocchi, Stirling and others, as well as the tangent
numbers, secant numbers, etc., are closely related to the
Bernoulli numbers. Following an idea due to J. Bernoulli it was
known that
$$S_n(k)=\frac{1}{n+1}\sum_{i=0}^n\binom{n+1}{i}B_ik^{n+1-i},
\text{ cf.[1, 2, 3, 4],} $$ where $B_i$ are Bernoulli numbers. Let
$n, k$ be positive integers ($k>1$), and let
$$T_n(k)=-1^n+2^n-3^n+4^n-5^n+\cdots+(-1)^{k-1}(k-1)^n.$$
In this paper we evaluate the alternating sums of powers of
consecutive integers. Following an idea due to Euler, we study a
formula for $T_n(k)$:
$$T_n(k)=\frac{(-1)^{k+1}}{2}\sum_{l=0}^{n-1}\binom nl E_lk^{n-l}+
\frac{E_n}{2}\left(1+(-1)^{k+1}\right),$$ where $E_l$ are the
Euler numbers.

\head 2. Sums of powers
\endhead

Euler numbers are defined by the generating function as follows:
$$G(t)=\frac{2}{e^t+1}=e^{Et}=\sum_{n=0}^{\infty}E_n \frac{t^n}{n!},\text{
 $|t|<\pi$, }\tag1$$ where we use the technique method notation by
 replacing $E^m$ by $E_m$ ($m\geq 0$), symbolically.
Let $x$ be the variable. Then we consider
$$F(t,x)=\frac{2}{e^t+1}e^{xt}=\sum_{n=0}^{\infty}E_n(x)\frac{t^n}{n!}.
\tag2$$ By (1) and (2), we see that
$$F(t,x)=G(t)e^{xt}=\sum_{n=0}^{\infty}\left(\sum_{k=0}^n\binom nk
E_k x^{n-k}\right) \frac{t^n}{n!}.$$ From this, we derive
$$E_n(x)=\sum_{n=0}^n\binom nk E_k x^{n-k}, $$
which is called the Euler polynomials. Note that $E_n(0)=E_n$.
From the definition of Euler numbers, we can derive  the below
relation:
$$ \left(E+1\right)^n+ E_n =2\delta_{0,n}, \text{ where $\delta_{0,n}$ is
Kronecker Symbol }. \tag3$$
 By (3), we easily see that
 $E_0=1$, $E_1=-\frac{1}{2},$ $E_2=0,$ $E_3=\frac{1}{4}, \cdots,$
 $E_{2k}=0.$
For any positive integer $n$, it is easy to show that
$$-2\sum_{l=0}^{\infty}(-1)^{l+n}e^{(l+n)t}+2\sum_{l=0}^{\infty}(-1)^le^{lt}
=2\sum_{l=0}^{n-1}(-1)^le^{lt}. \tag4$$ From (1), (2) and (4), we
derive
$$\sum_{m=0}^{\infty}\left(E_m+(-1)^{n+1}E_m(n)\right)\frac{t^m}{m!}
=\sum_{m=0}^{\infty}\left(2\sum_{l=0}^{n-1}(-1)^{l}l^n\right)\frac{t^m}{m!}.$$
Therefore we obtain the following theorem:
 \proclaim{ Theorem 1}
 Let $m, n$ be positive integers ($n>1$). Then we have
 $$\sum_{l=0}^{n-1}(-1)^ll^m=\frac{1}{2}\left((-1)^{n+1}E_m(n)+Em\right).$$
That is,
$$T_m(n)=\frac{(-1)^{n+1}}{2}\sum_{l=0}^{m-1}\binom ml E_l
n^{m-l}+\frac{E_m}{2}\left(1+(-1)^{n+1}\right).$$ If $n\equiv 0
(mod  2)$ then
$$T_m(n)+\frac{1}{2}\sum_{l=0}^{m-1}\binom ml E_l
n^{m-l}=0.$$
\endproclaim

\head 3. Euler-Zeta function
\endhead

Let $\Gamma(s)$ be the gamma function. Then we easily see that
$$\frac{1}{\Gamma(s)}\int_{0}^{\infty}t^{s-1}F(-t,x)dt=\sum_{n=0}^{\infty}\frac{(-1)^n2}{\Gamma(s)}
\int_{0}^{\infty}t^{s-1}e^{-(n+x)t}dt=2\sum_{n=0}^{\infty}\frac{(-1)^n}{(n+x)^s},
\text{ $s\in\Bbb C.$}$$ Thus, we can define the Euler-zeta
function follows:

\proclaim{Definition 2} For $s\in\Bbb C ,$ $x\in\Bbb R$ with
$0<x<1,$ define $$\zeta_E
(s,x)=2\sum_{n=0}^{\infty}\frac{(-1)^n}{(n+x)^s}, \text{ and, }
\zeta_E(s)=\zeta_{E}(s,1)=\sum_{n=1}^{\infty}\frac{(-1)^{n-1}}{n^s}.$$
\endproclaim

\proclaim{Lemma 3} Foe $s\in\Bbb C,$ we have
$$\zeta_{E}(s,x)=\frac{2}{\Gamma(s)}\int_{0}^{\infty}t^{s-1}\frac{e^{-xt}}{1+e^{-t}}dt.$$\endproclaim

From (2) and Lemma 3, we can derive the below Theorem 4.
\proclaim{ Theorem 4} Let $n$ be the positive integer. Then we
have $$\zeta_E (-n,x)=E_n(x).$$
\endproclaim
Remark. From Theorem 4, we note that
$$\aligned
&1=E_0(1)=\zeta_E(0)=2(1-1+1-1+\cdots ),\\
&\frac{1}{2}=E_1(1)=\zeta_E(-1)=2(1-2+3-4+\cdots),\\
&0=E_2(1)=\zeta_E(-2)=2(1^2-2^2+3^2-4^2+\cdots),\\
&\cdots.
\endaligned$$

\enddocument